# AN ITERATIVE PROCEDURE FOR GENERAL PROBABILITY MEASURES TO OBTAIN I-PROJECTIONS ONTO INTERSECTIONS OF CONVEX SETS

By Bhaskar Bhattacharya

*Southern Illinois University*

The iterative proportional fitting procedure (IPFP) was introduced formally by Deming and Stephan in 1940. For bivariate densities, this procedure has been investigated by Kullback and Rüschendorf. It is well known that the IPFP is a sequence of successive I-projections onto sets of probability measures with fixed marginals. However, when finding the I-projection onto the intersection of arbitrary closed, *convex* sets (e.g., marginal stochastic orders), a sequence of successive I-projections onto these sets may *not* lead to the actual solution. Addressing this situation, we present a new iterative I-projection algorithm. Under reasonable assumptions and using tools from Fenchel duality, convergence of this algorithm to the true solution is shown. The cases of infinite dimensional IPFP and marginal stochastic orders are worked out in this context.

**1. Introduction.** For two probability measures (PM) $P$ and $Q$ defined on an arbitrary measurable space $(\mathcal{X}, \mathcal{B})$, the *I-divergence* or the *Kullback–Leibler distance* between $P$ and $Q$ is defined as

$$I(P|Q) = \begin{cases} \int \ln(dP/dQ)\,dP, & \text{if } P \ll Q, \\ +\infty, & \text{otherwise.} \end{cases}$$

If $R$ is any PM with $P \ll R, Q \ll R$, then $I(P|Q)$ can equivalently be expressed as $I(P|Q) = \int (dP/dR) \ln((dP/dR)/(dQ/dR))\,dR$. Here and in the sequel we observe the conventions that $\ln 0 = -\infty, \ln(a/0) = +\infty, 0 \cdot (\pm\infty) = 0$.

Although $I(P|Q)$ is not a metric, it is always nonnegative and equals 0 if and only if $P = Q$ (a.e.). Hence, it is often interpreted as a measure

---

of "divergence" or "distance" between $P$ and $Q$. Other popular names of $I(P|Q)$ are information for discrimination, cross-entropy, information gain, and so on.

For a given $Q$ and a specified set of PM's $\mathcal{C}$, it is often of interest to find the $R \in \mathcal{C}$ which satisfies

$$(1.1) \qquad I(R|Q) = \inf_{P \in \mathcal{C}} \int \ln(dP/dQ)\, dP \ (< \infty).$$

This $R$ is called the *I-projection* of $Q$ onto $\mathcal{C}$. So, the I-projection, when defined, corresponds to a finite I-divergence. Csiszár [6] has shown that $R$ exists uniquely if $\mathcal{C}$ is variation-closed and there exists $P \in \mathcal{C}$ such that $I(P|Q) < \infty$. Csiszár [6] also gives a characterization of $R$ as follows: $R$ is the I-projection of $Q$ onto the convex set $\mathcal{C}$ if and only if

$$(1.2) \qquad I(P|Q) \geq I(P|R) + I(R|Q), \ \text{or}$$
$$\int \left(\frac{dP}{dQ} - \frac{dR}{dQ}\right) \ln\left(\frac{dR}{dQ}\right) dQ \geq 0,$$

for every $P \in \mathcal{C}$ (equality holds if $R$ is an algebraic inner point of $\mathcal{C}$).

If $Q$ is a finite measure but not a PM, with $c = \int dQ$ and $Q' = Q/c$, the I-divergence between a PM $P$ and $Q$ is given by $I(P|Q) = I(P|Q') - \ln c$. It can be seen easily that the same characterization result in (1.2) holds, that is, $Q$ and $Q'$ would lead to the same I-projection onto $\mathcal{C}$ although $I(P|Q)$ need not be nonnegative if $Q$ is not a PM.

I-projections play a key role in the information theoretic approach to statistics [12, 19]. The areas of maximization of entropy [16, 23] and the theory of large deviations [27] also use this concept. The iterative proportional fitting procedure (IPFP) [9, 18] commonly used in contingency tables is actually an iterative algorithm of successive I-projection problems [6, 7, 8, 15]. Maximum likelihood estimation in log-linear models for multinomial distributions is equivalent to solving an I-projection problem [11, 13, 14].

Depending on the form of the set $\mathcal{C}$, it may be difficult to find a solution to the I-projection problem in (1.1). In the discrete case, if $\mathcal{C}$ can be expressed as $\bigcap_{i=1}^{k} \mathcal{C}_i$ where each $\mathcal{C}_i$ is a closed, linear set [i.e., $P_1, P_2 \in \mathcal{C}_i \Rightarrow \alpha P_1 + (1-\alpha)P_2 \in \mathcal{C}_i$, for all $\alpha$ for which $\alpha P_1 + (1-\alpha)P_2$ is a PM], then Csiszár [6] has shown that the sequence of cyclic iterated I-projections onto individual $\mathcal{C}_i$ converges to the solution of (1.1). Dykstra [10] modified this procedure to work for the case when $\mathcal{C}_i$ are arbitrary closed, convex sets subject to a limiting condition which was later removed by Winkler [28]. Later Bhattacharya and Dykstra [3] interpreted Dykstra's procedure in the context of Fenchel duality.

In the infinite-dimensional bivariate case, the IPFP is fitting (adjusting) a PM to two given marginals. Kullback [20] considered this problem, and Rüschendorf [25] proved its convergence. When finding the I-projection of a general PM $Q$ on a convex set $\mathcal{C}$, Bhattacharya and Dykstra [2] used Fenchel duality to identify an equivalent dual problem which might be easier to solve (depending on $\mathcal{C}$) than the I-projection problem at hand. They worked out several examples to demonstrate the utility of the duality approach.

In this paper we consider the infinite-dimensional case when the constraint set can be expressed as an intersection of a finite number of arbitrary variation-closed, convex sets. The motivation for this paper comes from the fact that successive iterative I-projections, as considered by Rüschendorf [25], onto these sets may *not* lead to the actual solution (see Example 3.1 for an analytical example). The rest of the paper is organized as follows. In Section 2 we introduce the necessary notation and preliminaries to the problem. An iterative algorithm, a modified version of the discrete case [3, 10], is presented in Section 3 which is shown to converge to the correct solution. In dual formulation, a cyclic descent algorithm is obtained which amounts to minimizing over one variable at a time while all the others are held fixed. We establish the correspondence between the primal and dual solutions at every step of the two algorithms, which aids to prove the convergence of I-projections. In Section 4 we consider the case when marginal PM's are stochastically greater than or equal to the given PM's. The infinite-dimensional IPFP follows when each $\mathcal{C}_i$ is defined as a linear set of fixed marginals. The duality approach used in this paper to solve this problem seems simpler and more intuitive. This generalizes easily to the case when more than two variables are involved.

**2. Preliminaries.** We begin with some necessary notation. Let the underlying probability space be denoted by $(\Omega, \mathcal{F}, Q)$, where $\Omega$ is the sample space, $\mathcal{F}$ is the $\sigma$-field of subsets of $\Omega$, and $Q$ is a given PM defined on elements of $\mathcal{F}$. We will use the notation $\int f \, dQ$ to indicate $\int_\Omega f(\omega) \, dQ(\omega)$.

We will work with the normed, linear vector space $L_1(Q)$ since $I(P|Q) < \infty$ implies $dP/dQ \in L_1(Q)$. Consider the function $f$ defined as

$$(2.1) \qquad f(x) = \begin{cases} \int x \ln x \, dQ, & \text{if } x \geq 0, \int x \, dQ = 1, \\ +\infty, & \text{otherwise,} \end{cases}$$

and let

$$(2.2) \qquad \mathcal{C}_0 = \{x \in L_1(Q) : x = dP/dQ \text{ for some } P \in \mathcal{C}\}.$$

Then the (primal) problem in (1.1) can be expressed as

$$(2.3) \qquad \inf_{x \in \mathcal{C}_0} f(x).$$

Since the primal space is taken as $L_1(Q)$, the dual space is given by $L_\infty(Q)$ [22], which is the space of bounded functions. However, this space is too restrictive, as shown by Bhattacharya and Dykstra [2] (see also Remark 3.2 in Section 3). A more general dual space $\overline{M}(\Omega, \mathcal{F})$, the set of extended valued $\mathcal{F}$-measurable functions on $\Omega$, is more useful in this context. So $\overline{M}(\Omega, \mathcal{F})$ will be our dual space as well.

The *convex conjugate* of $f$, denoted by $f^*(y)$, for $y \in \overline{M}(\Omega, \mathcal{F})$, is defined as

$$(2.4) \qquad f^*(y) = \sup_{x \in dom(f)} \left[ \int xy \, dQ - f(x) \right],$$

where $dom(f) = \{x \in L_1(Q) : f(x) < \infty\}$. Bhattacharya and Dykstra [2] showed that, for the function $f$ given by (2.1), the convex conjugate is given by

$$(2.5) \qquad f^*(y) = \ln\left( \int e^y \, dQ \right)$$

($Q$ need not be a PM).

A subset $K$ of a vector space is said to be a *cone* if $x \in K \Rightarrow \alpha x \in K, \forall \alpha \geq 0$. For an arbitrary subset $S$ of $L_1(Q)$, the *conjugate* (or *dual*) cone of $S$ is given by

$$S^\oplus = \left\{ y \in \overline{M}(\Omega, \mathcal{F}) : \int xy \, dQ \geq 0, \text{ for all } x \in S \right\}.$$

The dual problem to (2.3) is given by $\inf_{y \in \mathcal{C}_0^\oplus} f^*(y)$. Theorems 2.1 and 2.2 [2] below present a sufficient condition which can be used to identify the solutions of the I-projection problem and its dual, and express one solution in terms of the other. Also, from (2.4) it follows that $f(x) + f^*(y) \geq \int xy \, dQ \geq 0, \ \forall x \in S, y \in S^\oplus$. Theorem 2.1 shows that if $f(x_0) + f^*(y_0) \leq 0$ for some $x_0 \in S, y_0 \in S^\oplus$, then $x_0, y_0$ solve the primal and dual optimization problems, respectively.

THEOREM 2.1.   *Assume $S$ is a subset of $L_1(Q)$ which intersects $dom(f)$ and $T \subset S^\oplus$. Then $x_0 \in S$ and $y_0 \in T$ are respective solutions of*

$$\inf_{x \in S} f(x) = \inf_{x \in S \cap dom(f)} \int x \ln x \, dQ \quad and \quad \inf_{y \in T} f^*(y) = \ln\left[ \inf_{y \in T} \int e^y \, dQ \right]$$

*if*

$$(2.6) \qquad f(x_0) + f^*(y_0) \leq 0.$$

*Moreover, in this situation $\inf_{x \in S} f(x) = -\inf_{y \in T} f^*(y)$.*

In the trivial case when $T = \{0\}$, it is easy to see that $y_0 = 0, f^*(y_0) = 0$. Then (2.6) gives that $f(x_0) \leq 0$. Since $f(x) \geq 0 \ \forall x$, we have $f(x_0) = 0$. From Theorem 2.1, $\inf_{x \in S} f(x) = -\inf_{y \in T} f^*(y) = 0$. Thus, we have $Q \in \mathcal{C}$.

In (2.5), substituting $y(\omega) = \alpha z(\omega), \alpha \geq 0, \omega \in \Omega$, for a fixed $z \in \overline{M}(\Omega, \mathcal{F})$, one may interpret $f^*(y)$ as the cumulant generating function of a random variable $z$ (with PM $Q$) when $f^*$ is evaluated along the one-dimensional path $y = \alpha z$. Thus, minimizing $f^*$ over the (one-dimensional) region $\mathcal{K}^{\oplus} = \{v : v(\omega) = \alpha z(\omega), \alpha \geq 0\}$ is equivalent to minimizing $I(P|Q)$ over the set $\mathcal{K} = \{x = dP/dQ \in L_1(Q) : \int zx\, dQ \geq 0\}$ since $\mathcal{K}^{\oplus}$ is the conjugate cone of $\mathcal{K}$. In other words, Theorem 2.1 shows that minimizing the cumulant generating function over nonnegative values of $\alpha$ is equivalent to finding the I-projection of $Q$ onto the set of distributions over $z(\omega)$ values with a nonnegative mean.

THEOREM 2.2. *Assume $S$ is a subset of $L_1(Q)$ and $y_0$ is a solution to* $\inf_{y \in S^{\oplus}} \int e^y\, dQ < \infty$. *Then $x_0 = e^{y_0} / \int e^{y_0}\, dQ$ is the solution to*

$$(2.7) \qquad \inf_{x \in S} f(x) = \inf_{x \in S \cap dom(f)} \int x \ln x\, dQ$$

*if $x_0 \in S$. If $S$ is either (1) convex, variation-closed and contained in $dom(f)$ or (2) a variation-closed, convex cone, then $x_0 \in S$, and, hence, $x_0$ must solve (2.7).*

Using Theorem 2.2, many exponential families of distributions are obtained as solutions to the I-projection problem when $S$ contains appropriate moment constraints [2]. From now on, to use the result of Theorem 2.2, we will replace $\mathcal{C}_0$ in (2.2) and (2.3) by

$$\mathcal{K}_0 = \{\alpha x : x \in \mathcal{C}_0, \alpha \geq 0\}.$$

**3. The algorithm and its convergence.** We assume that $Q$ is a given PM and there exists $V \in \mathcal{C}$ such that $I(V|Q) < \infty$. We wish to find the I-projection of $Q$ onto $\mathcal{C} = \bigcap_{i=1}^{t} \mathcal{C}_i$, where each $\mathcal{C}_i$ is a variation-closed, convex set of PM's. We first state the algorithm in the primal form, and then in the dual form. The main point is that successive iterative projections may *not* work for general convex sets (see Example 3.1 below). Hence, we propose some adjustments to be made before an I-projection is done.

*Primal formulation of the algorithm*:

- Initialization: Set $S_{0,i} = P_{0,i} = Q$ and begin with $n = 1, i = 1$.
- Implementation:

1. Set

$$(3.1) \qquad \begin{aligned} \frac{dS_{n,1}}{dQ} &= \frac{dP_{n-1,t}}{dQ} \left( \frac{dP_{n-1,1}}{dS_{n-1,1}} \right)^{-1}, \qquad \text{for } i = 1, \quad \text{and} \\ \frac{dS_{n,i}}{dQ} &= \frac{dP_{n,i-1}}{dQ} \left( \frac{dP_{n-1,i}}{dS_{n-1,i}} \right)^{-1}, \qquad \text{for } i = 2, \ldots, t \end{aligned}$$

(using the conventions described in the Introduction), and find $P_{n,i} = \pi_i(S_{n,i})$ [where $\pi_i(S_{n,i})$ is the I-projection of $S_{n,i}$ onto $\mathcal{C}_i$].

2. If $i < t$, increase $i$ by 1 and go to step 1 above. If $i = t$, increment $n$ by 1, set $i = 1$ and go to step 1 above.

In words, to implement the algorithm for the primal problem, at the first cycle we find the I-projection of $P_{1,i-1}$ onto $\mathcal{C}_i$, namely, $P_{1,i}$ for $1 \le i \le t$ (with $P_{1,0} = Q$). At the $n$th cycle $i$th step $(n \ge 2, 1 \le i \le t)$, we first form $dS_{n,i}/dQ$ by adjusting the current I-projection density $dP_{n-1,t}/dQ$ for $i = 1$, $dP_{n,i-1}/dQ$ for $i \ge 2$. The adjustment amounts to taking out the effect of the last I-projection from the previous cycle, which is $dP_{n-1,i}/dS_{n-1,i}$. Thus, although at the first cycle $S_{1,1} = Q, S_{1,i} = P_{1,i-1}, i \ge 2$, are all PM's, for $n \ge 2$, due to the adjustment made, the measure $S_{n,i}$ need not be a PM. Since the algorithm finds the I-projection from $S_{n,i}$, we need $S_{n,i}$ to be a finite measure so that (1.2) may be used in the proof of convergence of the algorithm. Thus, the following assumption is made:

ASSUMPTION A.   $\sup_{n,i} \int dS_{n,i} = M_1 < \infty$.

We show below (following Lemma 3.1) that if each $\mathcal{C}_i$ is a *linear* set, then no adjustment is necessary, $S_{n,1} = P_{n-1,t}, S_{n,i} = P_{n,i-1}, i \ge 2, \ \forall n$ are all PM's, and Assumption A clearly holds. In general, however, it is not easy to identify situations when Assumption A may be violated. One can monitor the value of $M_1$ while executing the algorithm (see the discussion following Theorem 3.1). One faces similar difficulty in the discrete case [10]. In the remarks following Example 3.2, different choices of $\Omega$ and $\mathcal{C}_i$'s are considered to discuss Assumption A.

We show in Lemma 3.1 that the I-projections $P_{n,i}$ and the corresponding densities $dP_{n,i}/dS_{n,i}$ defined by the algorithm exist $\forall n, i$. In the discrete case, Dykstra [10] assumed that the individual I-projections exist. Later Winkler [28] proved such existence; however, his proof depends on the discreteness of the sample space.

LEMMA 3.1.   *We assume that $Q$ is a given PM and there exists $V \in \mathcal{C}$ such that $I(V|Q) < \infty$. Then the densities $dP_{n,i}/dS_{n,i}$ defined by the algorithm exist $\forall n, i$.*

PROOF.   We prove this by induction. First, we show that, in the first cycle, all I-projections exist. By assumption, the projection of $Q$ onto $\mathcal{C}_1$, that is, $P_{1,1}$, exists. Applying (1.2) with $P = V \in \mathcal{C}_1$ (since $V \in \bigcap_{i=1}^k \mathcal{C}_i$) and $R = P_{1,1}$, we get $I(V|P_{1,1}) \le I(V|Q) < \infty$. Since $V \in \mathcal{C}_2$ also, it follows that $P_{1,2}$ exists. Thus, continuing in this way, it follows that each of $P_{1,3}, \dots, P_{1,k}$ exists. Now assume that all I-projections exist up to the $n$th cycle, $(i-1)$st step. To show that $P_{n,i}$ exists, it is enough to verify that $I(V|S_{n,i}) < \infty$ (since $V \in C_i$). It follows from (3.1) that $dS_{n,i}/dQ = \prod_{j=1(j \neq i)}^{t}(dP_{a,j}/dS_{a,j})$, where $a = n$ if $j \leq i$, otherwise $a = n - 1$. Then

$$I(V|S_{n,i}) = \int \ln\left(\frac{dV}{dS_{n,i}}\right) dV$$

$$= \int \ln\left(\frac{dV}{dQ}\right) dV - \int \ln\left(\frac{dS_{n,i}}{dQ}\right) dV$$

$$= I(V|Q) - \sum_{j=1(j \neq i)}^{t} \int \ln\left(\frac{dP_{a,j}}{dS_{a,j}}\right) dV$$

$$= I(V|Q) - \sum_{j=1(j \neq i)}^{t} \int \left(\frac{dV}{dQ} - \frac{dP_{a,j}}{dQ}\right) \ln\left(\frac{dP_{a,j}}{dS_{a,j}}\right) dQ$$

$$- \sum_{j=1(j \neq i)}^{t} I(P_{a,j}|S_{a,j})$$

$$\leq I(V|Q) - \sum_{j=1(j \neq i)}^{t} I(P_{a,j}|S_{a,j}),$$

since, by (1.2), every term in the middle sum is nonnegative (with $V = P, R = P_{a,j}, Q = S_{a,j}$) and using the induction hypothesis. Now replacing $S_{n,i}$ by $S'_{n,i} = S_{n,i}/\int dS_{n,i}$, we get

$$I(V|S'_{n,i}) + \sum_{j=1(j \neq i)}^{t} I(P_{a,j}|S'_{a,j}) \leq I(V|Q) + \sum_{j=1}^{t} \ln \int dS_{a,j}$$

$$\leq I(V|Q) + t \ln M_1 < \infty$$

by Assumption A. Thus, $I(V|S'_{n,i}) < \infty$, and consequently, $I(V|S_{n,i}) = I(V|S'_{n,i}) - \ln \int dS_{n,i} < \infty$ by Assumption A. Hence, $P_{n,i}$ exists, $I(P_{n,i}|S_{n,i}) < \infty$ and the densities $dP_{n,i}/dS_{n,i}$ exist $\forall n, i$. $\square$

One may also observe the following facts concerning the I-projections obtained from the algorithm. Using (1.2), it follows that, for $A \subset \Omega$, $\pi_i(Q)(A) = 0$ would imply either $Q(A) = 0$ or $P(A) = 0$ for all $P \in C_i$ [with $I(P|Q) < \infty$]. By reasoning inductively, it follows that $S_{n,i}(A) = 0$ would imply either $Q(A) = 0$ or $P(A) = 0$ for all $P \in C_i$ [with $I(P|S_{n,i}) < \infty$]. It also follows that once an I-projection assigns mass 0 (a.e. $Q$) to a set $A$ in $\Omega$, all subsequent I-projections also assign mass 0 (a.e. $Q$) to the same set $A$.

If the $C_i$'s were actually linear sets so that equality holds in (1.2), our procedure would reduce to the successive iterative projections algorithm and no adjustment is needed. To see this, let $P \in \mathcal{C}_i$. Then using (3.1),

$$
\begin{aligned}
I(P|S_{n,i}) &= \int \ln\left(\frac{dP}{dS_{n,i}}\right) dP \\
(3.2) \qquad &= \int \frac{dP}{dQ} \ln\left(\frac{dP/dQ}{dS_{n,i}/dQ}\right) dQ \\
&= \int \frac{dP}{dQ} \ln\left(\frac{dP/dQ}{dP_{n,i-1}/dQ}\right) dQ + \int \frac{dP}{dQ} \ln\left(\frac{dP_{n-1,i}}{dS_{n-1,i}}\right) dQ.
\end{aligned}
$$

From (1.2), setting $Q = S_{n-1,i}$, $R = P_{n-1,i}$, it follows that

$$
\int \left(\frac{dP}{dS_{n-1,i}} - \frac{dP_{n-1,i}}{dS_{n-1,i}}\right) \ln\left(\frac{dP_{n-1,i}}{dS_{n-1,i}}\right) dS_{n-1,i} = 0,
$$

or $\int (dP/dQ - dP_{n-1,i}/dQ) \ln(dP_{n-1,i}/dS_{n-1,i}) \, dQ = 0$. Hence, the last line of (3.2) is equal to

$$
\int \frac{dP}{dQ} \ln\left(\frac{dP/dQ}{dP_{n,i-1}/dQ}\right) dQ + \int \frac{dP_{n-1,i}}{dQ} \ln\left(\frac{dP_{n-1,i}}{dS_{n-1,i}}\right) dQ
$$

$$
= I(P|P_{n,i-1}) + I(P_{n-1,i}|S_{n-1,i}).
$$

Clearly, the $P \in \mathcal{C}_i$ which minimizes $I(P|S_{n,i})$ is the same as the $P$ which minimizes $I(P|P_{n,i-1})$. Thus, here we can take $S_{n,i} = P_{n,i-1}$, and, consequently, Assumption A is always true.

The following simple example demonstrates that the successive iterative projections procedure does not work for general convex sets.

EXAMPLE 3.1.   Suppose we like to find the I-projection of $Q$ [= the uniform distribution on $(0,1)$] onto the class of all distributions on $(0,1)$ whose first two moments are at least 0.7. The constraint region can be expressed as $\mathcal{C}_1 \cap \mathcal{C}_2$, where $\mathcal{C}_1 = \{P : E_P(X) \geq 0.7\}$, $\mathcal{C}_2 = \{P : E_P(X^2) \geq 0.7\}$, where $X$ is a random variable on $(0,1)$. It can be easily seen that the successive iterative projections algorithm produces the solution $R$ with density $r(x) = e^{2.672x + 1.943x^2}/17.120$, $0 < x < 1$. However, the correct solution is $R^*$ with density $r^*(x) = e^{3.932x^2}/7.845$, $0 < x < 1$.

In the dual formulation, the algorithm is a cyclic descent algorithm which successively minimizes over one function at a time while the others are held fixed. The dual problem is equivalent to

$$
\inf_{y \in (\bigcap_{i=1}^{t} \mathcal{K}_i)^{\oplus}} \int e^y \, dQ = \inf_{y \in \mathrm{cl}(\mathcal{K}_1^{\oplus} + \cdots + \mathcal{K}_t^{\oplus})} \int e^y \, dQ
$$

$$
= \inf_{y_i \in \mathcal{K}_i^{\oplus}, 1 \leq i \leq t} \int e^{y_1 + \cdots + y_t} \, dQ.
$$

Typically the dual constraint region would be the direct sum $(\mathcal{K}_1^{\oplus} + \cdots + \mathcal{K}_t^{\oplus})$. However, the direct sum of closed sets need not be closed [1], and hence our dual constraint region is the closure of the direct sum. At each step of the primal algorithm, a dual problem can be identified using Theorem 2.2. Following the duality approach, not only does the process seems simpler and more intuitive, but also, depending on the constraints, the dual problem may be more tractable than the I-projection problem.

*Dual formulation of the algorithm*:

- Initialization: Set $y_{0,i} = 0$ and begin with $n = 1, i = 1$.
- Implementation:

  1. Let $y_{n,i}$ denote the solution to

  $$\inf_{y \in \mathcal{K}_i^{\oplus}} \int e^{y_{n,1} + \cdots + y_{n,i-1} + y + y_{n-1,i+1} + \cdots + y_{n-1,t}} \, dQ.$$

  2. If $i < t$, increase $i$ by 1 and go to step 1. If $i = t$, increment $n$ by 1, set $i = 1$ and go to step 1.

In the following we present five lemmas which are crucial in proving the main result of this paper (Theorem 3.1). Recall, for a closed convex cone $\mathcal{K}$, the dual problem can be stated as

$$(3.3) \qquad \inf_{y \in \mathcal{K}^{\oplus}} \ln \int e^y \, dQ.$$

To begin, in Lemma 3.2 we establish a necessary condition for the solution of (3.3).

LEMMA 3.2. *Let $\mathcal{K}^{\oplus}$ be the conjugate of a closed convex cone $\mathcal{K}$. If $y_0$ solves the dual problem in (3.3), then*

$$(3.4) \qquad \begin{aligned} &\int y_0 e^{y_0} \, dQ = 0, \\ &\int y e^{y_0} \, dQ \geq 0 \qquad \forall y \in \mathcal{K}^{\oplus}. \end{aligned}$$

PROOF. Since $y_0$ minimizes $g(y) = \ln \int e^y \, dQ$ on the convex set $\mathcal{K}^{\oplus}$ and $g$ is Gateaux differentiable at $y_0$, we get $(d/d\alpha)g(\alpha y + (1-\alpha)y_0)|_{\alpha=0} \geq 0$. Applying this for the function $g(y)$, we get $(d/d\alpha) \ln \int e^{\alpha y + (1-\alpha)y_0} \, dQ|_{\alpha=0} \geq 0$, or

$$\int (y - y_0) \frac{e^{y_0}}{\int e^{y_0} \, dQ} \, dQ \geq 0 \qquad \forall y \in \mathcal{K}^{\oplus}.$$

By choosing $y = cy_0$ first with $c > 1$ and then with $c < 1$ (since $\mathcal{K}^{\oplus}$ is a cone), we obtain the results given in (3.4). □

Now we are able to relate the solutions from the primal and dual algorithms.

LEMMA 3.3.   *At the $(n,i)$th step, the solutions from the two algorithms are related by*

$$(3.5) \qquad \frac{dP_{n,i}}{dS_{n,i}} = \frac{e^{y_{n,i}}}{\int e^{y_{n,i}} \, dS_{n,i}} \quad and \quad I(P_{n,i}|S_{n,i}) = -\ln \int e^{y_{n,i}} \, dS_{n,i}.$$

*Also, $S_{n,i}$ and $P_{n,i}$ can be expressed as*

$$\frac{dS_{n,i}}{dQ} = c_{n,i} e^{\sum_{j=1}^{i-1} y_{n,j} + \sum_{j=i+1}^{t} y_{n-1,j}},$$

$$(3.6)$$

$$\frac{dP_{n,i}}{dQ} = \frac{e^{\sum_{j=1}^{i} y_{n,j} + \sum_{j=i+1}^{t} y_{n-1,j}}}{\int e^{\sum_{j=1}^{i} y_{n,j} + \sum_{j=i+1}^{t} y_{n-1,j}} \, dQ},$$

*where $c_{n,i}$ is given by*

$$(3.7) \qquad c_{n,i} = \left[ \int e^{\sum_{j=1}^{i-1} y_{n,j} + \sum_{j=i}^{t} y_{n-1,j}} \, dQ \right]^{-1}$$

$$\times \prod_{k=1}^{n-1} \left[ \frac{\int e^{\sum_{j=1}^{i} y_{k,j} + \sum_{j=i+1}^{t} y_{k-1,j}} \, dQ}{\int e^{\sum_{j=1}^{i-1} y_{k,j} + \sum_{j=i}^{t} y_{k-1,j}} \, dQ} \right],$$

$1 \le i \le t$ *(assuming $\sum_{j=1}^{0} y_{n,j} = \sum_{j=t+1}^{t} y_{n,j} = 0, y_{0,j} = 0$). In addition,*

$$(3.8) \qquad \int y_{n,i} \, dP_{n,i} = 0, \qquad \int y \, dP_{n,i} \ge 0 \qquad \forall \, y \in \mathcal{K}_i^{\oplus}.$$

*The I-divergences can be expressed in terms of dual solutions as*

$$I(P_{n,i}|S_{n,i}) = \sum_{k=1}^{n} -\ln \frac{\int e^{\sum_{j=1}^{i} y_{k,j} + \sum_{j=i+1}^{t} y_{k-1,j}} \, dQ}{\int e^{\sum_{j=1}^{i-1} y_{k,j} + \sum_{j=i}^{t} y_{k-1,j}} \, dQ},$$

$$(3.9)$$

$$\sum_{i=1}^{t} I(P_{n,i}|S_{n,i}) = -\ln \int e^{\sum_{j=1}^{t} y_{n,j}} \, dQ.$$

PROOF.   From Theorem 2.2, at the $(n,i)$th step of the primal algorithm, the solution to $\inf_{x \in \mathcal{K}_i} \int x \ln x \, dS_{n,i}$ is given by $dP_{n,i}/dS_{n,i} = e^{y_{n,i}} / \int e^{y_{n,i}} \, dS_{n,i}$, where $y_{n,i}$ solves the dual problem $\inf_{y \in \mathcal{K}_i^{\oplus}} \int e^y \, dS_{n,i}$. Also, by Theorem 2.1 we obtain $I(P_{n,i}|S_{n,i}) = -\ln \int e^{y_{n,i}} \, dS_{n,i}$.

We prove (3.6) and (3.7) by induction. For $n = 1$, we have $S_{1,i} = P_{1,i-1}$, and, consequently,

$$\frac{dP_{1,i}}{dS_{1,i}} = \frac{dP_{1,i}}{dP_{1,i-1}} = \frac{e^{y_{1,i}}}{\int e^{y_{1,i}} \, dP_{1,i-1}}.$$

Applying this equation recursively for $i-1, \ldots, 1$ and combining the resulting equations, we obtain $dP_{1,i}/dQ = e^{\sum_{j=1}^{i} y_{1,j}} / \int e^{\sum_{j=1}^{i} y_{1,j}} \, dQ$ with $c_{1,i} = (\int e^{\sum_{j=1}^{i-1} y_{1,j}} \, dQ)^{-1}$ for $1 \leq i \leq t$. Thus, (3.6) and (3.7) hold for $n=1$. Suppose (3.6) and (3.7) also hold up to step $(n, i-1)$ of the algorithm. To show these equations hold for the $(n, i)$th step, we start with the definition of $dS_{n,i}/dQ$ in (3.1) and then apply (3.5) for $(n-1, i)$ to obtain

$$
\begin{aligned}
\frac{dS_{n,i}}{dQ} &= \frac{dP_{n,i-1}}{dQ} \left( \frac{dP_{n-1,i}}{dS_{n-1,i}} \right)^{-1} \\
&= \frac{e^{\sum_{j=1}^{i-1} y_{n,j} + \sum_{j=i}^{t} y_{n-1,j}}}{\int e^{\sum_{j=1}^{i-1} y_{n,j} + \sum_{j=i}^{t} y_{n-1,j}} \, dQ} \frac{\int e^{y_{n-1,i}} \, dS_{n-1,i}}{e^{y_{n-1,i}}} \\
&= \frac{e^{\sum_{j=1}^{i-1} y_{n,j} + \sum_{j=i+1}^{t} y_{n-1,j}}}{\int e^{\sum_{j=1}^{i-1} y_{n,j} + \sum_{j=i}^{t} y_{n-1,j}} \, dQ} \\
&\quad \times \int e^{y_{n-1,i}} c_{n-1,i} e^{\sum_{j=1}^{i-1} y_{n-1,j} + \sum_{j=i+1}^{t} y_{n-2,j}} \, dQ \\
&= \frac{c_{n-1,i} \int e^{\sum_{j=1}^{i} y_{n-1,j} + \sum_{j=i+1}^{t} y_{n-2,j}} \, dQ}{\int e^{\sum_{j=1}^{i-1} y_{n,j} + \sum_{j=i}^{t} y_{n-1,j}} \, dQ} e^{\sum_{j=1}^{i-1} y_{n,j} + \sum_{j=i}^{t} y_{n-1,j}} \\
&= c_{n,i} e^{\sum_{j=1}^{i-1} y_{n,j} + \sum_{j=i}^{t} y_{n-1,j}},
\end{aligned}
$$

where $c_{n,i}$ is defined in (3.7). Also, using the expression for $dS_{n,i}/dQ$ derived above, we get

$$
\begin{aligned}
\frac{dP_{n,i}}{dQ} &= \frac{dP_{n,i}}{dS_{n,i}} \frac{dS_{n,i}}{dQ} \\
&= \frac{e^{y_{n,i}}}{\int e^{y_{n,i}} (dS_{n,i}/dQ) \, dQ} c_{n,i} e^{\sum_{j=1}^{i-1} y_{n,j} + \sum_{j=i}^{t} y_{n-1,j}} \\
&= \frac{e^{\sum_{j=1}^{i} y_{n,j} + \sum_{j=i+1}^{t} y_{n-1,j}}}{\int e^{\sum_{j=1}^{i} y_{n,j} + \sum_{j=i+1}^{t} y_{n-1,j}} \, dQ},
\end{aligned}
$$

which proves the desired results.

In (3.4), setting $Q = S_{n,i}, y_0 = y_{n,i}$, it follows that $\int y_{n,i} e^{y_{n,i}} \, dS_{n,i} = 0$, which, using (3.5), gives $\int y_{n,i} \, dP_{n,i} = 0$ as in (3.8). The rest of (3.8) follows from the fact that $dP_{n,i}/dQ \in \mathcal{K}_i$ and $y \in \mathcal{K}_i^{\oplus}$.

To derive (3.9), note that

$$
I(P_{n,i} | S_{n,i}) = -\ln \int e^{y_{n,i}} \, dS_{n,i} = -\ln \left( c_{n,i} \int e^{\sum_{j=1}^{i} y_{n,j} + \sum_{j=i+1}^{t} y_{n-1,j}} \, dQ \right).
$$

The expression for $I(P_{n,i}|S_{n,i})$, in (3.9), is obtained using the value of $c_{n,i}$ from (3.7). Also, writing $a_{k,i} = \int e^{\sum_{j=1}^{i} y_{k,j} + \sum_{j=i+1}^{t} y_{k-1,j}} \, dQ$, we get $\sum_{i=1}^{t} I(P_{n,i}|S_{n,i}) = \sum_{i=1}^{t} \sum_{k=1}^{n} -\ln(a_{k,i}/a_{k,i-1}) = -\ln a_{n,t}$ (assuming $a_{k,0} = a_{k-1,t}, a_{0,t} = 1$). This proves the lemma.  □

The next lemma essentially proves that the successive I-divergences are nondecreasing in nature.

LEMMA 3.4.  *Under Assumption* A:

(i)  $I(P_{n,i}|S_{n,i}) - I(P_{n-1,i}|S_{n-1,i}) \geq I(P_{n,i}|P_{n-1,i})$.
(ii)  $I(P_{n,i}|S_{n,i})$ *is nondecreasing in* $n$, *for every* $i$.

PROOF.  To prove (i), we may write

$$I(P_{n,i}|S_{n,i}) - I(P_{n-1,i}|S_{n-1,i})$$
$$= \int \ln\left(\frac{dP_{n,i}}{dS_{n,i}}\right) dP_{n,i} - \int \ln\left(\frac{dP_{n-1,i}}{dS_{n-1,i}}\right) dP_{n-1,i}$$
$$= \int \ln\left(\frac{dP_{n,i}}{dQ}\right) dP_{n,i} - \int \ln\left(\frac{dS_{n,i}}{dQ}\right) dP_{n,i} - \int \ln\left(\frac{dP_{n-1,i}}{dS_{n-1,i}}\right) dP_{n-1,i}$$
$$= \int \ln\left(\frac{dP_{n,i}}{dQ}\right) dP_{n,i} - \left\{ \int \ln\left(\frac{dP_{n,i-1}}{dQ}\right) dP_{n,i} - \int \ln\left(\frac{dP_{n-1,i}}{dS_{n-1,i}}\right) dP_{n,i} \right\}$$
$$\quad - \int \ln\left(\frac{dP_{n-1,i}}{dS_{n-1,i}}\right) dP_{n-1,i}$$
$$= \int \ln\left(\frac{dP_{n,i}}{dP_{n,i-1}}\right) dP_{n,i} + \int \left(\frac{dP_{n,i}}{dQ} - \frac{dP_{n-1,i}}{dQ}\right) \ln\left(\frac{dP_{n-1,i}}{dS_{n-1,i}}\right) dQ$$
$$\geq I(P_{n,i}|P_{n,i-1}),$$

since the second term is nonnegative by using the characterization of I-projection in (1.2) [with $P = P_{n,i}, R = P_{n-1,i}$ and $Q = S_{n-1,i}$ in the first line of (1.2), we can write $\int \ln(dP_{n,i}/dS_{n-1,i})(dP_{n,i}/dQ) \, dQ \geq \int \ln(dP_{n,i}/dP_{n-1,i})(dP_{n,i}/dQ) \, dQ + \int \ln(dP_{n-1,i}/dS_{n-1,i})(dP_{n-1,i}/dQ) \, dQ$, and the result follows by moving all terms to the left of the inequality and algebra] and using Assumption A. The case when $i = 1$ follows similarly.

Proof of (ii) follows from (i) since $I(P_{n,i}|P_{n,i-1}) \geq 0$.  □

The next lemma shows that the I-divergences obtained from the algorithm are also uniformly bounded above.

LEMMA 3.5.  *Under Assumption* A, $I(P_{n,i}|S_{n,i})$ *are bounded above uniformly in* $n, i$.

Proof. By assumption, $\exists V \in \mathcal{C}$ such that $I(V|Q) < \infty$. Writing $dP_{n,i}/dQ = (dP_{n,i}/dS_{n,i})(dS_{n,i}/dQ)$ and using (3.1), it follows that $dP_{n,i}/dQ = \prod_{j=1}^{t}(dP_{a,j}/dS_{a,j})$, where $a = n$ if $j \leq i$, and otherwise $a = n - 1$. Then

$$
\begin{aligned}
I(V|P_{n,i}) &= \int \ln\left(\frac{dV}{dP_{n,i}}\right) dV \\
&= \int \ln\left(\frac{dV}{dQ}\right) dV - \int \ln\left(\frac{dP_{n,i}}{dQ}\right) dV \\
&= I(V|Q) - \sum_{j=1}^{t} \int \ln\left(\frac{dP_{a,j}}{dS_{a,j}}\right) dV \\
&= I(V|Q) - \sum_{j=1}^{t} \int \left(\frac{dV}{dQ} - \frac{dP_{a,j}}{dQ}\right) \ln\left(\frac{dP_{a,j}}{dS_{a,j}}\right) dQ \\
&\quad - \sum_{j=1}^{t} I(P_{a,j}|S_{a,j}) \\
&\leq I(V|Q) - \sum_{j=1}^{t} I(P_{a,j}|S_{a,j}),
\end{aligned}
$$

(3.10)

since by (1.2) [with $P = V, R = P_{a,j}$ and $Q = S_{a,j}$ in (1.2) and decomposing (1.2) as shown in the proof of Lemma 3.4] every term in the middle sum is nonnegative under Assumption A. Thus,

$$
(3.11) \qquad I(V|P_{n,i}) + \sum_{j=1}^{t} I(P_{a,j}|S_{a,j}) \leq I(V|Q) < \infty.
$$

Thus, $I(P_{n,i}|S_{n,i})$ are bounded above uniformly in $n, i$.  □

Thus, Lemmas 3.4 and 3.5 together imply that $I(P_{n,i}|P_{n,i-1}) \to 0$ as $n \to \infty$. The following assumption is needed in the proofs of Lemma 3.6 and Theorem 3.1.

ASSUMPTION B.   $\sup_{n,i} \int (\ln \frac{dS_{n,i}}{dQ}) dP_{n,i} = M_2 < \infty$.

Since $I(P_{n,i}|S_{n,i}) = I(P_{n,i}|Q) - \int \ln(dS_{n,i}/dQ) dP_{n,i}$, Assumption B would hold if one could show $I(P_{n,i}|Q)$ are uniformly bounded above. In the case of two linear set constraints of fixed marginals, Rüschendorf ([25], Lemma 4.4) shows this to be true. But for general convex sets, this does not follow easily and it is also difficult to identify situations when $M_2 = \infty$ may occur. Of course, one can monitor the value of $M_2$ along with $M_1$ while executing the algorithm (see the discussion following Theorem 3.1). In the Remarks 3.1–3.7 we discuss Assumption B in several cases involving $\Omega$ and $\mathcal{C}_i$.

The next lemma establishes the uniform integrability of the sequence of functions $e^{\sum_{j=1}^{t} y_{a,j}}$ (where $a = n$ if $j \le i$, otherwise $a = n - 1$) obtained from the dual algorithm.

LEMMA 3.6. *Under Assumption* B, *the sequence of functions* $e^{\sum_{j=1}^{t} y_{a,j}}$ *is uniformly integrable where* $a = n$ *if* $j \le i$, *and otherwise* $a = n - 1$.

PROOF. From (3.9) we have

$$I(P_{n,i}|S_{n,i}) - I(P_{n-1,i}|S_{n-1,i})$$
$$= -\ln \int e^{\sum_{j=1}^{i} y_{n,j} + \sum_{j=i+1}^{t} y_{n-1,j}} \, dQ + \ln \int e^{\sum_{j=1}^{i-1} y_{n,j} + \sum_{j=i}^{t} y_{n-1,j}} \, dQ,$$

which is nonnegative by Lemma 3.4(i), or

$$(3.12) \quad \int e^{\sum_{j=1}^{i} y_{n,j} + \sum_{j=i+1}^{t} y_{n-1,j}} \, dQ \le \int e^{\sum_{j=1}^{i-1} y_{n,j} + \sum_{j=i}^{t} y_{n-1,j}} \, dQ.$$

Thus, using (3.12), each of the terms in square brackets under the $\prod$ sign in (3.7) is less than (or equal to) 1 for $k = 1, \ldots, n - 1$. Since

$$\left( \int e^{\sum_{j=1}^{i-1} y_{n,j} + \sum_{j=i}^{t} y_{n-1,j}} \, dQ \right)^{-1} \le \left( \int e^{\sum_{j=1}^{t} y_{n,j}} \, dQ \right)^{-1}$$
$$= \exp\left( \sum_{i=1}^{t} I(P_{n,i}|S_{n,i}) \right)$$

by (3.12) and (3.9), the leading term of $c_{n,i}$ in (3.7) is also uniformly bounded above from Lemma 3.5. Thus, we conclude

$$(3.13) \qquad\qquad \sup_{n,i} c_{n,i} < \infty.$$

Using (3.6) and (3.8), we can write

$$\int \sum_{j=1}^{t} y_{a,j} e^{\sum_{j=1}^{t} y_{a,j}} \, dQ$$
$$= \left( \int \sum_{j=1}^{t} y_{a,j} \, dP_{n,i} \right) \left( \int e^{\sum_{j=1}^{t} y_{a,j}} \, dQ \right)$$
$$= \left( \int \left( y_{n,i} + \ln \frac{dS_{n,i}}{dQ} - \ln c_{n,i} \right) dP_{n,i} \right) \left( \int e^{\sum_{j=1}^{t} y_{a,j}} \, dQ \right)$$
$$= \left( \int \left( \ln \frac{dS_{n,i}}{dQ} - \ln c_{n,i} \right) dP_{n,i} \right) \left( \int e^{\sum_{j=1}^{t} y_{a,j}} \, dQ \right).$$

Now for the continuous convex function $\phi(x) = x \ln x$, we have $\lim_{x \to \infty}(\phi(x)/x) = \infty$. Then using a modified form of the criteria of Valle Poussain [21], it is enough to establish that

$$\sup_{n,i} \int \phi\big(e^{\sum_{j=1}^{t} y_{a,j}}\big)\, dQ < \infty.$$

However,

$$\sup_{n,i} \int \phi\big(e^{\sum_{j=1}^{t} y_{a,j}}\big)\, dQ$$

$$= \sup_{n,i} \int \sum_{j=1}^{t} y_{a,j} e^{\sum_{j=1}^{t} y_{a,j}}\, dQ$$

$$= \sup_{n,i} \Big( \int \Big( \ln \frac{dS_{n,i}}{dQ} - \ln c_{n,i} \Big)\, dP_{n,i} \int e^{\sum_{j=1}^{t} y_{a,j}}\, dQ \Big)$$

$$\leq \Big( \sup_{n,i} \int \Big( \ln \frac{dS_{n,i}}{dQ} \Big)\, dP_{n,i} + \sup_{n,i} \ln c_{n,i} \Big)$$

$$\times \sup_{n,i} \Big( \int e^{\sum_{j=1}^{t} y_{a,j}}\, dQ \Big),$$

which is finite by Assumption B, (3.13) and the fact that $\int e^{\sum_{j=1}^{t} y_{a,j}}\, dQ$ is nonincreasing in $i$ (and $n$) from (3.12). Hence, the result follows. $\quad\square$

The main focus is on the following theorem, which proves convergence of the I-projection solutions obtained from the algorithm described above by using the connection between the primal and dual solutions at every step. Since we are considering in (2.3) the infimum of a convex function over a closed convex region, the solution exists and is unique (see also [6]).

THEOREM 3.1.    *We assume that $Q$ is a given PM and there exists $V \in \mathcal{C}$ such that $I(V|Q) < \infty$. We also assume that the sequence of primal solutions $dP_{n,i}/dQ \in L_1(Q)$, Assumptions* A *and* B *hold and the corresponding dual solutions $y_{n,i} \in \overline{M}(\Omega, \mathcal{F})$ and $y_{n,i}$ are uniformly integrable. Then as $n \to \infty$, we have the following:*

1. *There exist (unique) $x_0 \in L_1(Q)$, $y_0 \in \overline{M}(\Omega, \mathcal{F})$ such that $\int x_0\, dQ$ is a PM, $dP_{n,i}/dQ \xrightarrow{w} x_0$; $x_0$, $y_0$ solve the primal and dual problems, respectively, and $x_0 = e^{y_0}/\int e^{y_0}\, dQ$.*
2. *Moreover, if $x_0 = dP^*/dQ$ $(P^* \in \mathcal{C})$, then $I(P^*|P_{n,i}) \to 0 \ \forall\, i$. Also,*

$$\|P_{n,i} - P^*\| = \int \Big| \frac{dP_{n,i}}{dQ} - \frac{dP^*}{dQ} \Big|\, dQ \to 0 \qquad \forall\, i,$$

   *where $\|P_{n,i} - P^*\|$ is the total variation distance between $P_{n,i}$ and $P^*$.*

3. $\sum_{i=1}^{t} I(P_{n,i}|S_{n,i}) \to I(P^*|Q)$.

PROOF. 1. Using Lemmas 3.4 and 3.5, $I(P_{n,i}|P_{n,i-1}) \to 0$ as $n \to \infty$, hence, by using the well-known relation [5, 6]

$$(3.14) \qquad \|P - R\| \le (2I(P|R))^{1/2},$$

for any two PM's $P$ and $R$, it follows that $\|P_{n,i} - P_{n,i-1}\| \to 0$ for all $2 \le i \le t$, that is,

$$(3.15) \qquad \int \left| \frac{dP_{n,i}}{dQ} - \frac{dP_{n,i-1}}{dQ} \right| dQ \to 0.$$

Using (3.8), we have

$$
\begin{aligned}
(3.16) \qquad \int \left( \sum_{i=1}^{t} y_{n,i} \right) \frac{dP_{n,t}}{dQ}\, dQ &= \sum_{i=1}^{t} \int y_{n,i} \left( \frac{dP_{n,t}}{dQ} - \frac{dP_{n,i}}{dQ} \right) dQ \\
&= \sum_{i=1}^{t} \int y_{n,i} \left\{ \sum_{s=i+1}^{t} \left( \frac{dP_{n,s}}{dQ} - \frac{dP_{n,s-1}}{dQ} \right) \right\} dQ.
\end{aligned}
$$

From (3.15), $|dP_{n,t}/dQ - dP_{n,i}/dQ| \to 0$ (a.e. $Q$) and it follows that $dP_{n,t}/dQ \in \bigcap_{i=1}^{t} \mathcal{K}_i$ (a.e. $Q$) for sufficiently large $n$. Since $\sum_{i=1}^{t} y_{n,i} \in \bigoplus_{i=1}^{t} \mathcal{K}_i^{\oplus}$, we get

$$\liminf_{n \to \infty} \int \left( \sum_{i=1}^{t} y_{n,i} \right) \frac{dP_{n,t}}{dQ}\, dQ \ge 0.$$

The (a.e. $Q$) boundedness of $|dP_{n,t}/dQ - dP_{n,i}/dQ|$ and uniform integrability of $y_{n,i}$ imply from (3.16) that

$$\limsup_{n \to \infty} \int \left( \sum_{i=1}^{t} y_{n,i} \right) \frac{dP_{n,t}}{dQ}\, dQ = \limsup_{n \to \infty} \int \left( \sum_{i=1}^{t} y_{n,i} \right)(0)\, dQ = 0.$$

Hence,

$$\lim_{n \to \infty} \int \left( \sum_{i=1}^{t} y_{n,i} \right) \frac{dP_{n,t}}{dQ}\, dQ = 0.$$

Adding and subtracting $\ln \int e^{\sum_{i=1}^{t} y_{n,i}}\, dQ$ to the left-hand side of (3.16), we get

$$(3.17) \qquad \lim_{n \to \infty} \left( f \left( \frac{dP_{n,t}}{dQ} \right) + f^* \left( \sum_{i=1}^{t} y_{n,i} \right) \right) = 0.$$

It follows from Lemma 3.6 that $e^{\sum_{i=1}^{t} y_{n,i}}$ is uniformly integrable. Thus, $e^{\sum_{i=1}^{t} y_{n,i}}$ is tight, and hence, relatively compact ([4], pages 35–41). So, given any sequence of increasing positive numbers, there is a subsequence, say, $\{n_j\}$, and an element $x_0 \in \overline{M}(\Omega, \mathcal{F})$ such that

$$(3.18) \qquad e^{\sum_{i=1}^{t} y_{n_j, i}} \xrightarrow{w^*} x_0$$

[1, 22]. This means

$$(3.19) \qquad \int e^{\sum_{i=1}^{t} y_{n_j, i}} x \, dQ \to \int x_0 x \, dQ \qquad \forall x \in L_1(Q)$$

as $j \to \infty$. Also, $x \in L_\infty(Q) \Rightarrow x \in L_1(Q)$. Using the convergence result in (3.19) only for $x \in \overline{M}(\Omega, \mathcal{F})$, we get the weak convergence of $e^{\sum_{i=1}^{t} y_{n_j, i}}$ in $L_1(Q)$. Thus

$$(3.20) \qquad e^{\sum_{i=1}^{t} y_{n_j, i}} \xrightarrow{w} x_0 \qquad \text{in } L_1(Q).$$

Using the constant function $1 \in \overline{M}(\Omega, \mathcal{F})$, we obtain

$$(3.21) \qquad \int e^{\sum_{i=1}^{t} y_{n_j, i}} \, dQ \to \int x_0 \, dQ$$

as $j \to \infty$. From Jensen's inequality and uniform integrability of $y_{n, i}$, it follows that

$$\int e^{\sum_{i=1}^{t} y_{n_j, i}} \, dQ \geq e^{\int \sum_{i=1}^{t} y_{n_j, i} \, dQ} > 0.$$

Hence, from (3.21), $\int x_0 \, dQ > 0$.

Combining (3.20) and (3.21), we obtain

$$(3.22) \qquad \frac{dP_{n_j, t}}{dQ} = \frac{e^{\sum_{i=1}^{t} y_{n_j, i}}}{\int e^{\sum_{i=1}^{t} y_{n_j, i}} \, dQ} \xrightarrow{w} \frac{x_0}{\int x_0 \, dQ}.$$

Since $f$ is lower semicontinuous, we have

$$(3.23) \qquad \liminf_{n_j \to \infty} f\left(\frac{dP_{n_j, t}}{dQ}\right) \geq f\left(\frac{x_0}{\int x_0 \, dQ}\right).$$

Since $f$ is convex and $\mathcal{K}$ is a closed convex cone, the solution to $\inf_{x \in \mathcal{K}} f(x)$ exists uniquely. From (2.6), if we can show for some $x_0 \in \mathcal{K}, y_0 \in \mathcal{K}^\oplus$, $f(x_0/\int x_0 \, dQ) + f^*(y_0) \leq 0$, then $x_0/\int x_0 \, dQ, y_0$ are (unique) solutions to the primal and dual problems, respectively. We show below that $x_0/\int x_0 \, dQ$ in (3.22) can be used for this purpose. So if the sequence $dP_{n, t}/dQ$ has many limit points, they must be equal (a.e. $Q$). Let $y_0 = \ln x_0$. Then using (3.17), (3.21) and (3.23), one gets

$$f\left(\frac{x_0}{\int x_0 \, dQ}\right) + f^*(y_0) = f\left(\frac{x_0}{\int x_0 \, dQ}\right) + \ln \int x_0 \, dQ$$

$$\leq \liminf_{n_j \to \infty} f\left(\frac{dP_{n_j,t}}{dQ}\right) + \lim_{n \to \infty} \ln \int e^{\sum_{i=1}^t y_{n_j,i}}\, dQ$$

(3.24)

$$\leq \liminf_{n_j \to \infty}\left\{ f\left(\frac{dP_{n_j,t}}{dQ}\right) + f^*\left(\sum_{i=1}^t y_{n_j,i}\right)\right\}$$

$$= 0.$$

To complete the proof, we need to show that $x_0 \in \bigcap_{i=1}^t \mathcal{K}_i$ and $y_0 \in (\bigcap_{i=1}^t \mathcal{K}_i)^\oplus$.

Using (3.6) and (3.8) for $i = t$, we get $\int y\, dP_{n_j,t} \geq 0$ for $y \in \mathcal{K}_t^\oplus$, or

(3.25) $$\int e^{\sum_{i=1}^t y_{n_j,i}} y\, dQ \geq 0 \qquad \forall y \in \mathcal{K}_t^\oplus.$$

Also, by (3.20)

(3.26) $$\int e^{\sum_{i=1}^t y_{n_j,i}} y\, dQ \to \int x_0 y\, dQ \qquad \forall y \in \mathcal{K}_t^\oplus$$

as $j \to \infty$. Combining (3.25) and (3.26), we obtain $\int x_0 y\, dQ \geq 0$, $\forall y \in \mathcal{K}_t^\oplus$, or equivalently, $x_0 \in \mathcal{K}_t$.

For any $1 \leq i \leq t-1$, with $y \in \mathcal{K}_i^\oplus$, by using (3.8), we can say

$$0 \leq \int y \frac{dP_{n_j,i}}{dQ}\, dQ = \int y \sum_{a=i}^{t-1}\left\{\frac{dP_{n_j,a}}{dQ} - \frac{dP_{n_j,a+1}}{dQ}\right\} dQ + \int y \frac{dP_{n_j,t}}{dQ}\, dQ.$$

Since

$$\limsup_{j \to \infty} \int y \sum_{a=i}^{t-1}\left\{\frac{dP_{n_j,a}}{dQ} - \frac{dP_{n_j,a+1}}{dQ}\right\} dQ = \limsup_{j \to \infty} \int y(0)\, dQ = 0$$

and by using (3.22),

$$0 \leq \int y \frac{dP_{n_j,t}}{dQ}\, dQ \to \int y \frac{x_0}{\int x_0\, dQ}\, dQ$$

(as $j \to \infty$) for all $y \in \mathcal{K}_i^\oplus$. Hence, $x_0 \in \mathcal{K}_i$, $1 \leq i \leq t-1$. Hence, $x_0 \in \bigcap_{i=1}^t \mathcal{K}_i$.

Now we show that $y_0 \in (\bigcap_{i=1}^t \mathcal{K}_i)^\oplus$. From (3.18), it follows that for the subsequence $\{n_j\}$ there exists $\tilde{y}$ such that

(3.27) $$\sum_{i=1}^t y_{n_j,i} \xrightarrow{w^*} \tilde{y}.$$

By lower semicontinuity of $f^*$,

(3.28) $$\liminf_{n_j \to \infty} f^*\left(\sum_{i=1}^t y_{n_j,i}\right) \geq f^*(\tilde{y}).$$

Using (3.28) and arguments similar to (3.24), it follows that

$$f\left(\frac{x_0}{\int x_0 \, dQ}\right) + f^*(\tilde{y}) \leq 0.$$

For $x \in \bigcap_{i=1}^t \mathcal{K}_i$, it follows that $\int (\sum_{i=1}^t y_{n_j,i}) x \, dQ \geq 0$ since $\sum_{i=1}^t y_{n_j,i} \in (\bigcap_{i=1}^t \mathcal{K}_i)^{\oplus}$. By (3.27), $\int x\tilde{y} \geq 0$. Hence, $\tilde{y} \in (\bigcap_{i=1}^t \mathcal{K}_i)^{\oplus}$, and $x_0/\int x_0 \, dQ$, $\tilde{y}$ solve the primal, dual problems, respectively. Since we have already established that $x_0/\int x_0 \, dQ, y_0$ solve the primal, dual problems, respectively, we must have $y_0 = \tilde{y}$ (a.e. $Q$).

2. For $P^* \in \mathcal{C}$, if $dP^*/dQ = x_0$ (a.e. $Q$), then

$$\begin{aligned}
I(P^*|P_{n,i}) &= \int \ln\left(\frac{dP^*}{dP_{n,i}}\right) dP^* \\
&= \int \left(\ln\frac{dP^*}{dQ} - \ln\frac{dP_{n,i}}{dQ}\right) dP^* \\
&\to \int \left(\ln\frac{dP^*}{dQ} - \ln\frac{dP^*}{dQ}\right) dP^* = 0
\end{aligned}$$

(as $n \to \infty$) by continuity of the ln function. By (3.14), it follows that $P_{n,i}$ converges to $P^*$ in total variation.

3. From (3.10), we have

$$I(V|P_{n,i}) = I(V|Q) - \sum_{j=1}^t \int \left(\frac{dV}{dQ} - \frac{dP_{a,j}}{dQ}\right) \ln\left(\frac{dP_{a,j}}{dS_{a,j}}\right) dQ - \sum_{j=1}^t I(P_{a,j}|S_{a,j}).$$

Setting $V = P^*$, since $dP_{a,j}/dQ \xrightarrow{w} dP^*/dQ$ and $I(P^*|P_{n,i}) \to 0$ as $n \to \infty$, we get the desired result. □

As mentioned earlier, we recommend that one can monitor the values of $M_1$ and $M_2$ while executing the algorithm. This can be done by inserting an extra step of computation while writing a computer program for the algorithm. If the algorithm is not going to converge correctly, then $M_1$ and/or $M_2$ should become excessively large. Otherwise the algorithm must converge to the correct solution. We find it extremely difficult to construct examples where either $M_1$ or $M_2$ might be infinite. Thus, it seems that such cases will only be possible under quite unusual circumstances (although we cannot prove or disprove this). We discuss the Assumptions A and B following Example 3.2 for different choices of $\Omega$ and $\mathcal{C}_i$'s in Remarks 3.1–3.4.

If $y_{n,i}$ were uniformly bounded [which would be true if we used $L_\infty(Q)$ as the dual space instead of $\overline{M}(\Omega, \mathcal{F})$], then the uniform integrability condition of $y_{n,i}$ is easily satisfied. Although the assumption that the $y_{n,i}$'s are uniformly integrable instead of uniformly bounded makes the proof of Theorem 3.1 rather complicated, this covers many situations which would not be possible otherwise, as the next example shows.

EXAMPLE 3.2. Suppose we wish to find the I-projection of the PM $Q$ with density (with respect to Lebesgue measure on the unit square) given by $q(x, y) = (4/5)(1 + xy), 0 < x, y < 1$, onto $\mathcal{C}_1 \cap \mathcal{C}_2$, where $\mathcal{C}_1 = \{P : E_P(\ln X) \geq -0.5\}, \mathcal{C}_2 = \{P : E_P(X + Y) \geq 1.3\}$. Writing $\mathcal{C}_1 = \{P : \int (0.5 + \ln x)(dP/dQ) dQ \geq 0\}$ and $\mathcal{C}_2 = \{P : \int (x + y - 1.3)(dP/dQ) dQ \geq 0\}$, it follows that the conjugate cones are given by $\mathcal{C}_1^{\oplus} = \{y(u, v) : y(u, v) = \alpha(0.5 + \ln u), \alpha \geq 0\}, \mathcal{C}_2^{\oplus} = \{y(u, v) : y(u, v) = \beta(u + v - 1.3), \beta \geq 0\}$. While checking the assumptions, we note that the $y_{n,1} (\in \mathcal{C}_1^{\oplus})$'s are unbounded, but uniformly integrable and the $y_{n,2} (\in \mathcal{C}_2^{\oplus})$'s are uniformly bounded (so uniformly integrable). After six cycles we obtain $\int dS_{6,1} = 0.8598, \int dS_{6,2} = 0.9049$ and $\int \ln(dS_{6,1}/dQ) dP_{6,1} = 2.81 \times 10^{-4}, \int \ln(dS_{6,2}/dQ) dP_{6,2} = -3.25 \times 10^{-4}$ [these integrals settle around these values for slightly higher $(n, i)$'s also]. Thus, we assume that Assumptions A and B hold. Applying the algorithm, after six cycles we obtain $P_{6,2}$, where $dP_{6,2}/d\lambda = e^{1.0394(u+v)} u^{0.3757}(1 + uv)/3.3451$. Noting $E_{P_{6,2}}(\ln X) = -0.4992$ and $E_{P_{6,2}}(X + Y) = 1.300$ and assuming that convergence is essentially obtained, one may take $P_{6,2}$ as a good approximation for $P^*$, the solution to the I-projection problem. In most cases (including this example) numerical integration techniques are needed for the above calculations.

It may be noted that many of the maximum entropy characterizations of families of distributions (e.g., [17]) can be obtained following our procedure. In Example 3.1, the solution $R^*$ have maximum entropy in the class of all distributions on $(0, 1)$ which have the first two moments at least 0.7.

REMARK 3.1. If the space $\Omega$ were discrete and finite, and the $\mathcal{C}_i$'s were linear, variation-closed sets, then our procedure would reduce to that of Csiszár [6] with the interpretation that $s_{n,1}(k) = p_{n-1,t}(k), s_{n,i}(k) = p_{n,i-1}(k)$, $n \geq 2$. Replacing the integrals by sums, here Assumption A holds easily since $\sum_k s_{n,1}(k) = 1, \forall n, i$. Also, $\sum_k \ln(s_{n,i}(k)/q(k))p_{n,i}(k) = \sum_k (\ln s_{n,i}(k))p_{n,i}(k) + \sum_k (-\ln q(k)) \times p_{n,i}(k) \leq \sum_k (-\ln q(k))p_{n,i}(k) \leq \sum_k (-\ln q(k)) < \infty$; thus Assumption B holds for $q(k) > 0$ (note $q(k) = 0$ implies $s_{n,i}(k) = p_{n,i}(k) = 0$, $\forall n, i$ with $0/0 = 0$; see the discussion following Lemma 3.1 and [3]). The pointwise convergence $p_{n,i}(k) \to p^*(k) \ \forall k$ is attained without assuming uniform integrability of the $y_{n,i}$, as can be seen by mimicking the arguments of Bhattacharya and Dykstra [3] when applied for the linear sets.

REMARK 3.2. If the space $\Omega$ were discrete and finite, and the $\mathcal{C}_i$'s were convex, variation-closed sets, our procedure would reduce to that of Bhattacharya and Dykstra [3] (also [10]). Here the $s_{n,i}$'s are only positive measures, and hence one must make Assumption A (noted by Dykstra [10] also). However, Assumption A implies that $s_{n,i}(k) \leq M_1 \ \forall k$ [since $s_{n,i}(k) > 0, \forall n, i$]. Then $\sum_k \ln(s_{n,i}(k)/q(k))p_{n,i}(k) = \sum_k \ln s_{n,i}(k)p_{n,i}(k) +$ $\sum_k(-\ln q(k))p_{n,i}(k) \le \ln M_1 + \sum_k(-\ln q(k)) < \infty$ [for $q(k) > 0$, $\forall k$]; thus, Assumption B holds. One need not make the assumption of uniform integrability of $y_{n,i}$ as shown by Bhattacharya and Dykstra [3].

REMARK 3.3. If $\Omega$ were discrete but with infinite support (so Bhattacharya and Dykstra [3] does not apply), then Assumption A holds when the $\mathcal{C}_i$'s are linear and must be assumed when the $\mathcal{C}_i$'s are convex. However, in general, it seems intractable to show that $M_2 < \infty$. One of the possible reasons is that $q(k)$ could be positive but infinitely small for many $k$. Consequently, Assumption B has to be made. Also, uniform integrability of $y_{n,i}$ has to be assumed.

REMARK 3.4. If the space $\Omega$ were infinite-dimensional and the $\mathcal{C}_i$'s were linear, variation-closed sets, our procedure would generalize to the infinite-dimensional case of Csiszár [6] (which remains unsolved still now), with the understanding that $S_{n,1} = P_{n-1,t}, S_{n,i} = P_{n,i-1}, i \ge 2, \forall n$ (see the derivation following Lemma 3.1). In this case, Assumption A holds easily since $\int dS_{n,i} = 1, \forall n, i$. For checking Assumption B, first note that, for $i \ge 2$, $\int \ln(dP_{n,i-1}/dQ) \, dP_{n,i} = -I(P_{n,i}|P_{n,i-1}) + I(P_{n,i}|Q)$ (similarly for $i = 1$). Although $I(P_{n,i}|P_{n,i-1}) \to 0$ as $n \to \infty$, in general it seems to be intractable to show that $\sup_{n,i} I(P_{n,i}|Q) < \infty$. But using particular definitions of the sets $\mathcal{C}_i$, this may be possible; for example, Rüschendorf [25] imposes restrictions on marginal densities (see his conditions B2, B3 and Lemma 4.4) and verifies that $\sup_{n,i} I(P_{n,i}|Q) < \infty$ for IPFP. In general, the uniform integrability of $y_{n,i}$ also has to be assumed in this case. For IPFP, Rüschendorf [25] also made this assumption, but shows that this follows from his conditions B2, B3 on marginal densities (see his Proposition 3.2).

REMARK 3.5. The assumption that there exists some $V \in \mathcal{C}$ such that $I(V|Q) < \infty$, is the same as the assumption of Csiszár ([6] Theorem 2.1) that $S(Q, \infty) \cap \mathcal{C} \ne \varnothing$, where $S(Q, \rho) = \{P : I(P|Q) < \rho\}, 0 < \rho \le \infty$, is the I-sphere with center $Q$ and radius $\rho$. Without this assumption, one may not be able to find an individual I-projection [(1.1) and afterward] $P_{n,i}$, the terms $dS_{n,i}/dQ$ or $dP_{n,i}/dS_{n,i}$ may not be defined or be arbitrary, $M_1$ and/or $M_2$ may become excessively large, the algorithm may not converge, or any combination of these may occur. Also, if the assumption is not valid, then the final result (if found) will not correspond to an I-projection by (1.1).

REMARK 3.6. It is one of the difficulties of the iterative procedures in optimization methods to specify the distance of the current estimate from the actual solution since small changes in the objective function between successive steps of the iterative procedure do not guarantee that the actual optimal solution is near by. However, from Theorem 3.1, one can use the I-projection $P_{n,i}$ and the quantity $\sum_{i=1}^{t} I(P_{n,i}|S_{n,i})$ to estimate $P^*$ and $I(P^*|Q)$, respectively, to arbitrary accuracy by choosing $n$ sufficiently large. Yet another bound can be obtained in the following way. Let $\hat{P} \in \bigcap_{i=1}^{t} C_i$ be close (in I-divergence sense) to $P_{n,t}$, and $P^*$ be the actual solution to our problem. We also assume that $\hat{P}$ is *known* a priori, or may be obtained by least square projection onto a subset of $\bigcap_{i=1}^{t} C_i$ or by some other method. Since by substituting $V = P^*$ in (3.11), we get $I(P^*|Q) \geq \sum_{i=1}^{t} I(P_{n,i}|S_{n,i})$ for sufficiently large $n$, it follows that

$$I(\hat{P}|Q) - \sum_{i=1}^{t} I(P_{n,i}|S_{n,i}) \geq I(\hat{P}|Q) - I(P^*|Q)$$

$$\geq I(\hat{P}|P^*) \geq \frac{1}{2}\left(\sum_{i=1}^{t} \int \left|\frac{d\hat{P}}{dQ} - \frac{dP^*}{dQ}\right| dQ\right)^2.$$

Consequently, using the leftmost term, we can specify an upper bound on both the I-divergence distance and the variation distance between the PM $\hat{P}$ and the true solution $P^*$. For instance, in Example 3.2 if one chooses the PM $\hat{P}$, where $d\hat{P}/d\lambda = ce^{\alpha(0.5+\ln u)+\beta(u+v-1.3)}(4/5)(1+uv)$ ($c =$ normalizing constant) with $\alpha = 0.4$, $\beta = 1.04$, then it is easily verified that $\hat{P}$ is in $C_1 \cap C_2$, and $I(\hat{P}|Q) - \sum_{i=1}^{2} I(P_{6,i}|S_{6,i}) = 0.0064$.

REMARK 3.7. The rate of convergence of the algorithm depends largely on the nature of the $C_i$'s. If all the $C_i$'s are orthogonal with each other (in terms of I-divergence), a single pass through each constraint will suffice. On the other hand, if one constraint has a narrow angle (in terms of I-divergence) with another, then many cycles will be needed to achieve a desired level of convergence (see [13]).

## 4. Marginal stochastic orders.

Since the introduction of the IPFP by Deming and Stephan [9], it has been widely used in many different fields. For a given bivariate density function, Kullback [20] considered the problem of matching (approximating) its marginal density functions to given univariate densities using an iterative algorithm (similar to IPFP in the discrete case). Recently, Rüschendorf [25] proved the convergence of this iterative algorithm under some regularity conditions. The purpose of this section is to view the problem considered by Rüschendorf [25] from a duality perspective, and to extend it to the case of (separate) row and column marginal stochastic orders. The process generalizes naturally to higher dimensions with restrictions involving possibly more than one marginal at a time.

Let $Q$ be a fixed, bivariate PM on $\mathcal{R}^2$ with $X$-marginal $Q_X$ and $Y$-marginal $Q_Y$. Also, let $G_X$ and $G_Y$ be given univariate PM's on $\mathcal{R}$ which are absolutely continuous with respect to $Q_X$ and $Q_Y$, respectively. We consider the problem of finding the I-projection of $Q$ onto $\mathcal{S}(G_X, G_Y)$, the class of all bivariate PM's whose $X$-marginal ($Y$-marginal) is stochastically greater than or equal to $G_X$ ($G_Y$). Thus, $\mathcal{S}(G_X, G_Y) = \mathcal{S}_1 \cap \mathcal{S}_2$, where $\mathcal{S}_1 = \{P : \int (I_{A_x} - G_X(x))(dP/dQ)\,dQ \leq 0, \forall\, x \in \mathcal{R}\}$, $\mathcal{S}_2 = \{P : \int (I_{B_y} - G_Y(y))(dP/dQ)\,dQ \leq 0, \ \forall\, y \in \mathcal{R}\}$ with $I_{A_x} = \{(u, v) : -\infty < u < x, -\infty < v \leq \infty\}$, $I_{B_y} = \{(u, v) : -\infty < u \leq \infty, -\infty < v < y\}$, $G_X(x) = Pr_{G_X}(X \leq x)$, $G_Y(y) = Pr_{G_Y}(Y \leq y)$.

Scrutinizing the definition of the conjugate cone, we see that the corresponding dual problem can be expressed as $\inf_{\mathcal{S}_1^{\oplus} \oplus \mathcal{S}_2^{\oplus}} \int e^y \, dQ$, where $y = y(u, v) = y_1(u) + y_2(v) \in \mathcal{S}_1^{\oplus} \oplus \mathcal{S}_2^{\oplus}$ and

$$\mathcal{S}_1^{\oplus} = \left\{ y(u, v) : y(u, v) = y_1(u), y_1(u) \text{ is nondecreasing, } \int y_1(u)\,dG_X(u) = 0 \right\},$$

$$\mathcal{S}_2^{\oplus} = \left\{ y(u, v) : y(u, v) = y_2(v), y_2(v) \text{ is nondecreasing, } \int y_2(v)\,dG_Y(v) = 0 \right\}.$$

Let $E_Q(g|\mathcal{I})$ denote the least square projection (or isotonic regression) of $g$ onto $\mathcal{I}$, the cone of all nondecreasing functions on $\mathcal{L}^2(Q)$ ([24], Chapter 8) with weights $Q$. Following the algorithm of Section 3, we first find the I-projection of $Q$ onto $\mathcal{S}_1$. Defining

$$
\begin{aligned}
y_{1,1}(u, v) = {} & \ln\left( E_Q\left( \frac{dG_X}{dQ_X} \Big| \mathcal{I} \right) \right)(u) \\
& - \int \ln\left( E_Q\left( \frac{dG_X}{dQ_X} \Big| \mathcal{I} \right) \right)(t)\,dG_X(t),
\end{aligned}
$$

it can be seen that $y_{1,1}(u, v)(= y_{1,1}(u)) \in \mathcal{S}_1^{\oplus}$, and for any $y(u, v)(= y(u)) \in \mathcal{S}_1^{\oplus}$, we have $(d/d\alpha) \int e^{y_{1,1} + \alpha(y - y_{1,1})}\,dQ|_{\alpha=0} = \int (y(t) - y_{1,1}(t))\,dG_X(t) = 0$, for $0 < \alpha < 1$, which implies that $y_{1,1}$ solves the dual problem. From Theorem 2.2, $E_Q(dG_X/dQ_X|\mathcal{I}) = a_1^*$ (say) solves the primal I-projection problem. If $(dQ/d(Q_X \otimes Q_Y))(u, v) = h(u, v)$ and $P_{n,i}^*$ is the I-projection at the $n$th cycle onto $\mathcal{S}_i$, then we can express $dP_{1,1}^*/d(Q_X \otimes Q_Y) = (dP_{1,1}^*/dQ)(dQ/d(Q_X \otimes Q_Y)) = a_1^* h$. Next we find the I-projection of $P_{1,1}^*$ onto $\mathcal{S}_2$. Following the last projection, it is given by $dP_{1,2}^*/dP_{1,1}^* = E_{P_{1,1}^{*Y}}(dG_Y/dP_{1,1}^{*Y}|\mathcal{I}) = b_1^*$ (say), where $P_{1,1}^{*Y}$ is the $Y$-marginal of $P_{1,1}^*$. Hence, we can write

$$
\begin{aligned}
\frac{dP_{1,2}^*}{d(Q_X \otimes Q_Y)} &= \frac{dP_{1,2}^*}{dP_{1,1}^*} \frac{dP_{1,1}^*}{dQ} \frac{dQ}{d(Q_X \otimes Q_Y)} \\
&= E_{P_{1,1}^{*Y}}\left( \frac{dG_Y}{dP_{1,1}^{*Y}} \Big| \mathcal{I} \right) E_Q\left( \frac{dG_X}{dQ_X} \Big| \mathcal{I} \right) \frac{dQ}{d(Q_X \otimes Q_Y)} \\
&= b_1^* a_1^* h.
\end{aligned}
$$

To begin cycle 2, we first form $W_{2,1}(A) = \int_A b_1^*(v)\,dQ$, $A \subset \mathcal{R}^2$. The I-projection, $P_{2,1}^*$, of $W_{2,1}$ onto $\mathcal{S}_1$ is given by $(dP_{2,1}^*/d(Q_X \otimes Q_Y)) = (dP_{2,1}^*/dW_{2,1})(dW_{2,1}/dQ)(dQ/d(Q_X \otimes Q_Y)) = a_2^* b_1^* h$. Thus, in general, after the $n$th cycle, $dP_{n,1}^*/d(Q_X \otimes Q_Y) = a_n^* b_{n-1}^* h$, $dP_{n,2}^*/d(Q_X \otimes Q_Y) = a_n^* b_n^* h$, where $a_n^*(u) = E_{P_{n-1,2}^{*X}}(dG_X/dP_{n-1,2}^{*X}|\mathcal{I})(u)$, $b_n^*(v) = E_{P_{n,1}^{*Y}}(dG_Y/dP_{n,1}^{*Y}|\mathcal{I})(v)$. By Theorem 3.1, there exist $a^*, b^*$ such that $a_n^* b_n^* \to a^* b^*$ (in variation), and we get $dP^*/d(Q_X \otimes Q_Y)(u,v) = a^*(u)b^*(v)h(u,v)$.

To verify Assumption A, we note that $\sup_{n,i} \int dW_{n,i} < \infty$ simplifies to $\sup_n \{\int a_n^* \, dQ_X, \int b_n^* \, dQ_Y\} < \infty$. Since

$$\int \ln(dS_{n,1}/dQ)\,dP_{n,1} = \int (\ln b_{n-1}^*)a_n^* b_{n-1}^* \, dQ$$
$$= \left(\int b_{n-1}^* \ln b_{n-1}^* \, dQ_Y\right)\left(\int a_n^* \, dQ_X\right),$$

$$\int \ln(dS_{n,2}/dQ)\,dP_{n,2} = \int (\ln a_n^*)a_n^* b_n^* \, dQ$$
$$= \left(\int a_n^* \ln a_n^* \, dQ_X\right)\left(\int b_n^* \, dQ_Y\right),$$

for Assumption B it is enough to assume that $\sup_n \{\int a_n^* \ln a_n^* \, dQ_X, \int b_n^* \times \ln b_n^* \, dQ_Y\} < \infty$ along with Assumption A. One also has to assume that $\ln a_n^*$ and $\ln b_n^*$ are uniformly integrable (a.e. $Q$). The stochastic ordering problem considered above would reduce to that of Rüschendorf [25] and Rüschendorf and Thomsen [26] if there were equality in the definitions of $\mathcal{S}_1$ and $\mathcal{S}_2$. The nondecreasing restrictions in the conjugate cones would be replaced by equalities, and in the algorithm, $W_{n,1}, W_{n,2}$ taken to be the same as $P_{n-1,2}^*, P_{n,1}^*$, respectively, $\forall n \geq 2$. Here, the above assumptions follow easily from the conditions imposed by Rüschendorf [25] on marginal densities (using $a_n^* = a_n, b_n^* = b_n$ and using 2.7 and 2.8 of [25]).

A modified set of Schrödinger (in)equalities in this situation can be expressed as

$$\int_{-\infty}^x \int_{-\infty}^\infty a^*(u)b^*(v)h(u,v)\,dv\,du \geq \int_{-\infty}^x g_X(u)\,du \qquad \forall x \in \mathcal{R},$$

$$\int_{-\infty}^\infty \int_{-\infty}^y a^*(u)b^*(v)h(u,v)\,dv\,du \geq \int_{-\infty}^y g_Y(v)\,dv \qquad \forall y \in \mathcal{R},$$

where $g_X, g_Y$ are fixed marginal densities. However, we note that they are not easy to solve for $a^*, b^*$, and the method described in this section would be needed to obtain the solution.

The above process can be extended to more than two variables, and the restrictions may involve more than one marginal at a time. The functions $a^*$ and $b^*$ will not be functions of $u$ and $v$ alone in this case. Theorem 3.1 yields the following.

THEOREM 4.1. *Let $Q$ be a $p$-variate PM. Let $\mathbf{m}_i \subset \{1, \ldots, p\}, 1 \leq i \leq k$, be subsets of indices. Let $Q_{\mathbf{m}_i}$ be the marginal PM of $Q$ for the indices in $\mathbf{m}_i$, and $\mathcal{S}_i$ be the closed convex cone that $Q_{\mathbf{m}_i}$ is stochastically greater than or equal to $G_{\mathbf{m}_i}$ (fixed). The I-projection $P^*$ of $Q$ onto $\mathcal{S} = \bigcap_{i=1}^m \mathcal{S}_i$ uniquely exists, and there exist $a_i$ such that $(dP^*/dQ)(u_1, \ldots, u_p) = a_1(u_{\mathbf{m}_1})a_2(u_{\mathbf{m}_2})\ldots a_k(u_{\mathbf{m}_k})$, where $u_{\mathbf{m}_i}$ is a vector with coordinates from $u_1, \ldots, u_p$ for indices in $\mathbf{m}_i$.*

As final comments, in this paper we have investigated the theoretical aspects of the proposed algorithm along with some simple examples to demonstrate its use. The cases of infinite-dimensional IPFP and marginal stochastic orders are illustrated from a duality perspective. The algorithm assumes that we are able to find the I-projection onto the individual sets $\mathcal{C}_i$. Simple cases, such as the ones in Examples 3.1 and 3.2, can be solved by Maple, Mathematica or IMSL routines. In other cases the level of difficulty will depend on our ability to compute the I-projection on the $\mathcal{C}_i$'s. Of course, with the invention of fast computing techniques, many such difficult tasks are well within reach. Thus, more research regarding computing may be needed to implement the algorithm in those cases.

**Acknowledgments.** I am very grateful to Professor Richard Dykstra for motivating and helping me during my Ph.D. studies to conduct research in this area. I am also grateful to two referees for their comments which have helped me to improve the presentation of this paper.

DEPARTMENT OF MATHEMATICS
SOUTHERN ILLINOIS UNIVERSITY AT CARBONDALE
CARBONDALE, ILLINOIS 62901-4408
USA
E-MAIL: bhaskar@math.siu.edu